# COMPUTER-ASSISTED INDEPENDENT STUDY IN MULTIVARIATE CALCULUS

L. Descalço[1], Paula Carvalho[1], J.P. Cruz[1], Paula Oliveira[1], Dina Seabra[2]

[1]*Departamento de Matemática, Universidade de Aveiro (PORTUGAL)*
[2]*Escola Superior de Tecnologia e Gestão de Águeda, Universidade de Aveiro (PORTUGAL)*

## Abstract

Learning mathematics requires students to work in an independent way which is particularly challenging for such an abstract subject. Advancements in technology and, taking the student as the focus of his own learning, led to a change of paradigm in the design and development of educational contents.

In this paper we describe the first experience with an interactive feedback and assessment tool (Siacua), based on parameterized math exercises, and explain how we use it to motivate student independent study in a multivariate calculus environment.

We have defined an index about the subject, trying to make it consensual enough for being used in other courses about multivariate calculus. Then we have created a concept map, selected some existing parameterized true/false questions from PmatE project and classified them using our concept map, for being reused in our system. For complementing the course we have created about one hundred parameterized multiple choice question templates in system Megua and generated about one thousand instances for using in Siacua.

Results based on data collected by this tool and also based on an informal survey are presented. This first experience allows us to conclude our approach has an important impact on student motivation and contributes to the success on learning multivariate calculus.

Keywords: Independent Study, Technology, Computer, Supported Education, Mathematics.

## 1 INTRODUCTION

The learning process, even when seen as independent study or active learning, is not necessarily an individual process; personal behavioural and environmental factors can influence people in a bi-directional way, even reciprocally.

Learning is not also an instantaneous process, it can occur without the manifestation of an observable change in behaviour and so achieved knowledge and the demonstration of this knowledge are distinct processes. This means students can achieve knowledge but not show it until they are motivated to do so.

At the same time that learning is, to a great extent, a result of watching the behaviour of models, students are more active and effective learners if they fell self-efficacy and self-confident in their ability to perform academic tasks successfully.

Some aspects of social learning, like online learning, participation in groups or discussion forums, namely in social networks, collaborative learning supported by study materials dedicated for specific topics, have shown some impact on education nowadays. It is in this context of social learning, although supported and complemented by the most classic and traditional methods, that this first experience in the University of Aveiro in Calculus with Several Variables, was done, based on this system of active and independent study.

The stage where this process takes place focuses essentially on three environments: the dedicated Web application Siacua (see [1, 2]), a thematic group in a social network and the classical classroom.

The first contact with the topics of the program is done at classroom in a more or less classic environment. The teacher provides knowledge and study material, illustrating with some examples and applications. The active and effective learning is left to the student assisted by the other two platforms.

The Web application Siacua, provides an environment where each student can access and practice all topics of the course and get immediate feedback on his performance provided by progress bars, associated with several concepts. The true/false and multiple choice questions are designed to provide the student some perception of his knowledge level in the subject matter by answering them. On the other hand, this application can be used as an effective tool for learning since it also gives the possibility of not answering the question and to see an answer with appropriate detail instead.

This individual work can be supplemented with discussion, request for clarification of doubts which may have arisen or remained, in the thematic discussion group on Facebook, where the students and teachers are together. In general, the act of verbalization or the obligation to write their doubts, are enough to help the students in understanding the problems and often helps to overcome small obstacles that inhibit learning.

## 2   PARAMETERIZED QUESTIONS

Learning mathematics in general requires a lot of work and sometimes routine tasks. The system we use allows the teacher to produce a large number of exercises about the same subject with the same difficulty level that can evaluate similar competencies for a large number of students. This is done by constructing parameterized questions.

In the University of Aveiro we can use the questions from PmatE (see [3, 4]), a project with 25 years, containing now many hundreds of parameterized questions in several areas of knowledge. Since the questions are parameterized, which means there are parameters in the questions that are instantiated in runtime by the computer, we have in fact many thousands of different questions available for use. We have selected and classified 22 questions from this system to use in our course about calculus with several variables.

Another project for creating parameterized question is Megua; see [5, 6]. It includes a package for Sage Notebook, allowing us to create multiple choice questions, with detailed solution, and to send them immediately to Siacua. For complementing the course we have created 117 parameterized multiple choice question templates in system Megua and used them to generate 1088 instances for using in Siacua.

## 3   USING THE BAYESIAN SYSTEM

We refer to our system as interactive learning, in the sense that it provides feedback in real-time and allows the student to decide what to do.

For computing this feedback, a user model based on Bayesian Networks (BN) is used, as described in [7]. This system has been tested with simulated students (see [8, 9]) and real students (see [10, 11]) with positive conclusions about student knowledge diagnosis.

We present a simple example to illustrate the use of this model in our course, in particular making clear the information teachers have to specify when authoring questions.

Let us assume that our course is named Calculus (C) and it can be divided in Derivatives (D), Integrals (I) and Partial Differential Equations (E). The graph representing the course is illustrated in Figure 1. In general, the concept map of a course is a tree, where the root is the course itself.

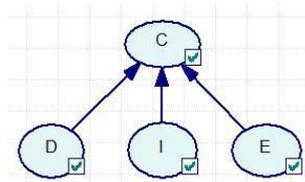

**Figure 1: Concept nodes**

The nodes in the graph correspond to binary random variables. To define the BN we have now to set the probabilities associated to each node. In our example we set

$$P(D=1)=0.5, P(D=0)=0.5, P(I=1)=0.5, P(I=0)=0.5, P(E=1)=0.5, P(E=0)=0.5$$

These probabilities represent the beliefs. For example, $P(D=1)=0.5$ means that we belief 50% that the student knows the subject Derivatives. In general, the leaves of the tree (nodes where no edges arrive) representing concepts are initiated with probability 0.5, meaning that without any evidence about the student knowledge, our beliefs are 50%.

We have then to set the conditional probabilities for the other nodes of the graph (non leaves). That is done automatically just by defining the importance, or weight, of each subtopic. In our example, assume that Derivatives has weight 0.5, Integrals weight 0.3 and PDEs weight 0.2. Hence the conditional probabilities are defined in Table 1, where, for example, $P(C=1|D=1, I=1, D=0)$, represents the belief about student knowledge in Calculus, assuming knowledge in Derivatives and Integrals but no knowledge in Partial Differential Equations. In general, the conditional probabilities are obtained by adding the corresponding weights.

| P(C=1 | | | | | | | | |
|---|---|---|---|---|---|---|---|---|
| D= | 1 | | | | 0 | | | |
| ,I= | 1 | | 0 | | 1 | | 0 | |
| ,E= | 1 | 0 | 1 | 0 | 1 | 0 | 1 | 0 |
| )= | 1.0 | 0.8 | 0.7 | 0.5 | 0.5 | 0.3 | 0.2 | 0.0 |

**Table 1: Conditional probabilities on concepts**

Suppose the student gives a correct answer for question $Q$ that has to do with Derivatives and also Integrals. Then a new node ($Q$) is added to the graph, as illustrated in Figure 2.

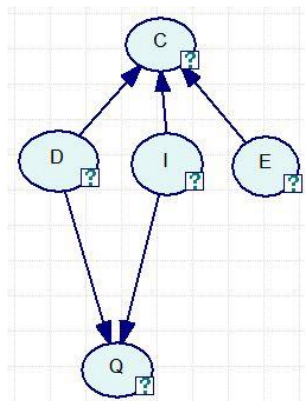

**Figure 2: Adding an evidence**

For adding this node to the Bayesian Network, the teacher has to specify the concepts it involves and also the corresponding weights. Question $Q$ may be, for example, 60% about derivatives and 40% about integrals. Moreover, the teacher has to specify some other parameters for computing the conditional probabilities in node $Q$. These parameters are, parameter "guess", the probability of answering correctly if the student does not know about Derivatives neither Integrals, which is normally 0.25 in multiple choice questions with 4 choices. Parameter "slip", representing the probability of giving a wrong answer knowing about both subjects Derivatives and Integrals. Also "level", the difficulty level (1-5) and "discr", the discrimination factor, most be provided.

Let us assume that, in our example guess=0.25, slip=0.2, level=1, discr=0.3. Using these parameters, the following probabilities are immediately set:

$$P(Q=1|D=0, I=0) = 0.25 \text{ (guess)} \quad P(Q=1|D=1, I=1) = 0.8 \text{ (1-slip)}$$

The remaining probabilities most be set, namely the probability of answering correctly (wrongly) knowing one concept but not the other. They are computed from these small set of parameters using the student model described in [7].

In practice, the teacher has only to provide the information described on the Megua code illustrated in Figure 3 before sending the question from Megua to Siacua.

```
SIACUAstart

level=1;  slip= 0.2; guess=0.25; discr = 0.3

concepts = [(D, 0.6), (I, 0.4)]

SIACUAend
```

**Figure 3: Parameters in Megua**

We observe that these parameters are specified by the teacher based on his/her teaching experience. But, by analysing the data being collected by Siacua application and comparing with the final marks students achieve in the course, we intend to improve the values of all question parameters.

We have defined an index with the relevant concepts for our course about calculus with several variables. We tried to do this choice as generic as possible in order to be consensual and useful for other similar courses. Hence, we have identified 57 "concepts" and defined the corresponding weights, obtaining a Bayesian network with graph illustrated by Figure 4.

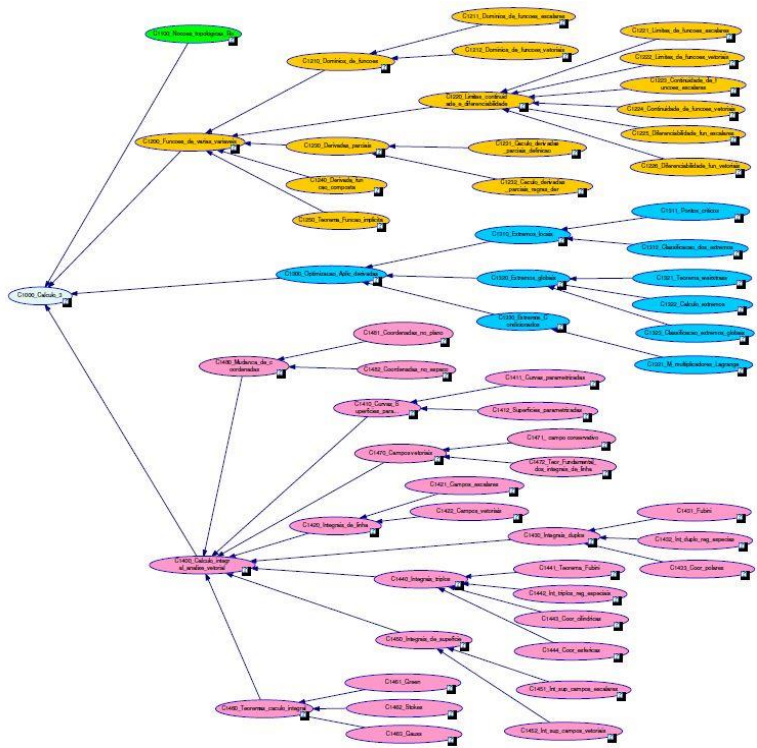

**Figure 4: Concepts in the calculus course**

The Bayesian system is used in a simple way. Each student is associated with a Bayesian network. Before any interaction, this network has only concepts. Then, each time the student answers a question, a new node is added to his/her network and the knowledge is propagated. This is done using the algorithms provided in Genie/Smile; see [12, 13]. Hence, students see several progress bars changing each time they answer a question, up to the top progress bar representing the course, which also changes a little after each answer.

For simplicity, if a student answers the same question a second time, we simply change the evidence in the already present node in the graph. So, the maximum number of nodes in the graph is limited to the number of concepts plus the number of available questions. Although all answers are kept by the system, only the last answer to each question is considered for diagnosis.

## 4   THE WEB APPLICATION

The Web application Siacua contains all exercises of our calculus course. Some of them are made available in PDF format with the corresponding solutions and are accessible through a simple side menu in the application. The interesting part, also accessible in the side menu, is the independent study area.

The independent study system is an open system in the sense that the student can see his progress, which is shown in the form of progress bars. The interaction is minimal, consisting of selecting questions and answering true/false or multiple choice questions as they appear. The selection of questions can be made by clicking the corresponding progress bar, what randomly selects a question related with the corresponding concept, or by introducing the question number. These question numbers are only shown in the moment the student is answering the question. This allows the student to identify the question for answering again latter but mainly for interacting with teachers and colleagues.

We aim to infer the cognitive state of the student, from the collected data, diagnosing his state in the learning process. Since the data is incomplete and sometimes not correct, this is a difficult step to implement in the interactive learning system. Nevertheless, this application provides a very simple way of presenting questions to the student, and giving some feedback. The student knows immediately if his answer is right or wrong, he can see the detailed solution and has information about his general progress in the course, given by the progress bars. Independent study system is illustrated in Figure 5.

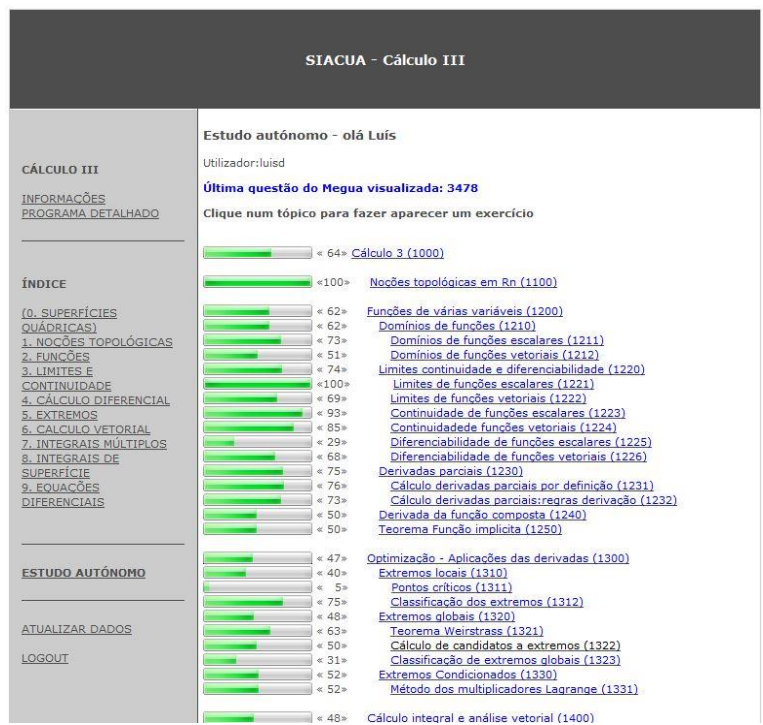

**Figure 5: Independent study**

The application can be used in computers and mobile devices. Students use it as a tool for studying and teachers can use it in the classes, as well. Moreover, a teacher can see the progress of each student in the application and also the average progress of the all class in each topic being studied, computed by the Bayesian networks.

## 5 FIRST EXPERIENCE WITH THE SYSTEM

### 5.1 Participants and procedure

As long as this system is being developed we applied it as an experiment in a course curricula of various engineering courses with a large number of students: about 500 students. The subject is Multivariate Calculus contemplating derivatives of functions of several variables, multiple integration and vector calculus. The students could use the Siacua application as a tool for independent study permanently but we have selected a topic to make possible the comparison of the students' performance in the exam. The selected topic was *local extrema* (classifying critical points of a function using the Hessian matrix). For a week students were induced to learn this topic using intensively the

Siacua system and a week later they were submitted to a classical exam with a problem (*test question*) of classifying the critical points of a given function.

## 5.2 Results

We observed that the number of students answering correctly the test question was very high when compared with other questions on the exam. Also the results obtained with the Web application on the selected topic were very high, as can be seen in Table 2.

| Number of students | Above 50% | Above 80% |
|---|---|---|
| 317 | 96% | 76% |

**Table 2: Success in the Web application on local extrema**

In Table 3 we summarize the use of Siacua by the whole population of 317 students that where submitted to assessment in Calculus. Column 1 shows the number of answers to multiple choice questions during the use of the application. We have an average of 45 answers per student but with a very high standard deviation, with some students answering many questions and some answering very few questions. Second column shows the answers to true/false questions where a similar situation occurs. Third and fourth columns show the number of logins in the system and the number of times the students enter the independent study area where questions with feedback are available.

|  | Megua answers | PmatE answers | Logins | Indep. study |
|---|---|---|---|---|
| Total | 15869 | 1912 | 14414 | 10845 |
| Average | 45,1 | 4,9 | 40,8 | 30,6 |
| Std. Dev. | 44 | 13 | 35,3 | 29,3 |

**Table 3: Student access on the platform**

Finally, Table 4 shows a comparison between the success in the question on the exam about local extrema and the mark in the whole exam, for a sample of 145 students randomly selected.

| Total number of students | Positive marks in test question | Positive marks in the exam |
|---|---|---|
| 145 | 81% | 64% |

**Table 4: Rate of success of test question comparison with all the others questions in the exam**

## 5.3 Analysis

With this experimental results we tend to believe that this tool is very efficient in the motivation and effective progress of student's independent study.

The students have been told that this is an ongoing project and application Siacua is a prototype. Nevertheless the feedback was very positive with many students thanking the teachers and incentivizing them to continue the development of exercises for Siacua.

Although we have not found any significant correlation between the student classification in Calculus and the knowledge computed by the Bayesian model, we conjecture, based on the very high standard deviations in Table 3, the main reason is that many students use the system just for reading the question's detailed solutions (so they use it as a tool for study) and so they try to guess the answers.

Nevertheless for 25,5% of population the difference between the classification in calculus and the Bayesian classification is less than 5% and for 46,8% of population is less than 10%.

Although these numbers have no precise meaning we guess some students use the system with care and for those there is some correlation between the two marks. We plan to use the system again as an optional assessment element with appropriate weight, so that the interested students use it more carefully.

In Table 5 we show the results of the informal survey, where 302 students have participated. First column contains the questions and second column contains students answers average in the usual Likert scale 1-5 where 1 stands for *completely disagree*, 3 is the neutral answer and 5 stands for *completely agree*.

| Q1 | I have used Siacua regularly to study Calculus 3 | 3,23 |
| Q2 | The true/false questions from PmatE are useful | 2,73 |
| Q3 | Multiple choice questions from Megua are useful | 3,21 |
| Q4 | It is better to have a dedicated application, like Siacua, for autonomous study, than to have the study material available on the internet (Moodle). | 4,46 |
| Q5 | The use of a platform like Siacua, contributes to a better learning. | 4,11 |
| Q6 | It is important that detailed solutions form Megua questions are available. | 4,39 |
| Q7 | It is important to know immediately if my answer to questions are right or wrong. | 4,54 |
| Q8 | It is important to have some information about my progress in the several topics studied. | 4,11 |
| Q9 | I think that this platform should not be used to assessment. | 2,36 |
| Q10 | Although there is Siacua for helping autonomous study, classes are still important. | 4,39 |
| Q11 | I would like to have more courses using this system. | 4,00 |
| Q12 | It is important that, in the future, the system makes suggestions of useful exercises. | 4,31 |

**Table 5: Informal survey**

From this survey we see that students value a dedicated Web application for helping independent study and in particular, the feedback provided by the system. The strongest agreement in the answers is for the importance of detailed solutions of exercises, then the importance of a dedicated Web application for independent study, and after that comes the importance of lectures. This informal survey, together with the feedback we receive form the students in classes, makes us believe our work has a positive impact on the teaching-learning process, and is an important complement to the classes, for independent study and motivation.

## 6   CONCLUSIONS

The results of this first experience, the feedback shown by the survey as well as the personal contact with students in classroom during the academic semester, shows that this project is appreciated by students and must be continued.

Application Siacua, still a prototype, can be improved in many ways. An obvious one is the inclusion of a mechanism for giving suggestions to the student based on the state of his interaction with the system, instead of simply showing him/her the state of his/her knowledge, that is, a system that continuously analyses each student's performance and is able to guide or direct him to the next step accordingly. This means that the student can progress to more difficult and complex tasks or, when a student fails in a particular set of questions associated with a specific concept, the system will be able to redirect him/her to the corresponding concept and also can deliver particular content in order to increase the student's proficiency.

This may include sending him to specific topics of another course, and so requires adding other layers of information to the concept map, like pre-requisite, proximity and dependence relations, other than just aggregation, what we have now. Hence also the student model must be improved in order to take into account the new information, including learning styles and many other relevant items in order for the system to suggest and even adapt appropriately.

Finally, by using the system for assessment in the future, we will collect valuable data generated by the student's interaction that will be used to adjust the Bayesian model parameters associated with the multiple choice questions and concepts. This will improve the concept weights and question classifications and, as a consequence, will improve the diagnosis in future uses.

## 7  AKNOWLEDGMENTS

This work was supported by CIDMA ("Center for Research & Development in Mathematics and Applications'") and FCT ("FCT- Fundação para a Ciência e a Tecnologia") through project UID/MAT/04106 /2013.

## REFERENCES

[1] Siacua: Interactive Computer Learning System, University of Aveiro (2015). http://siacua.web.ua.pt.

[2] Fonseca, M. M. G. (2014) Modelo Bayesiano do aluno no cálculo com várias variáveis. MSc Thesis. Dapartamento de Matemática - Universidade de Aveiro, Portugal.

[3] PmatE: Projeto Matemática Ensino (2015). http://pmate.ua.pt.

[4] Oliveira, M. Paula; Silva, Sabrina Vieira (2006). An overview of PmatE: developing software for all degrees of teaching. Proceedings of the International Conference in Mathematics Sciences and Sciences Education, June 11-14, University of Aveiro.

[5] Cruz, Pedro; Oliveira, Paula; Seabra, Dina (2012). Exercise templates with Sage. Tbilisi Mathematical Journal, 5(2), pp. 37-44.

[6] MEGUA package for parameterized exercises in Sage Mathematics (2015). http://cms.ua.pt/megua/.

[7] Descalço, L., Carvalho, Paula; Cruz, J.P.; Oliveira, Paula; Seabra, Dina. Using Bayesian networks and parameterized questions for independent study. Submitted to EDULEARN15.

[8] Millán, Eva (2000). Thesis: Bayesian system for student modeling. AI Commun. 13(4), pp. 277-278.

[9] Millán, Eva; Pérez-de-la-Cruz, José-Luis (2002). A Bayesian Diagnostic Algorithm for Student Modeling and its Evaluation. User Model. User-Adapt. Interact. 12(2-3), pp. 281-330.

[10] Castillo, Gladys; Descalço, L.; Diogo, Sandra; Millán, Eva; Oliveira, Paula; Anjo, Batel (2010). Computerized evaluation and diagnosis of student´s knowledge based on Bayesian Networks, Sustaining Tel: From Innovation to Learning and Practice, Lecture Notes in Computer Science, Volume 6383/2010, pp. 494-49.

[11] Castillo, Gladys; Millán, Eva; Descalço, L.; Oliveira, Paula; Diogo, Sandra (2013). Using bayesian networks to improve knowledge assessment. Computers & Education 60(1), pp. 436-447.

[12] Druzdzel, Marek J. (1999) SMILE: Structural Modeling, Inference, and Learning Engine and GeNIe: A development environment for graphical decision-theoretic models (Intelligent Systems Demonstration). In Proceedings of the Sixteenth National Conference on Artificial Intelligence (AAAI-99), pp. 902-903, AAAI Press/The MIT Press, Menlo Park, CA.

[13] Genie & Smile (2015), Decision Systems Laboratory of the University of Pittsburgh. https://dslpitt.org/genie/.